\documentclass[12pt]{article}
\usepackage{graphicx}
\usepackage{amsmath}
\usepackage{amsfonts}
\usepackage{amssymb}
\usepackage{epsfig,euscript}

\def \N{\mathbb{N}}

\def \R{\mathbb{R}}

\def \DR {Descartes' Rule of Signs}

\newtheorem{defn}{Definition}
\newtheorem{lem}{Lemma}
\newtheorem{prop}{Proposition}
\newtheorem{thm}{Theorem}

\newtheorem{notn}{Notation}
\newtheorem{rem}{Remark}

\newenvironment{proof}[1]{
  \trivlist \item[\hskip \labelsep{\it #1}]}{\hfill\mbox{$\square$}
  \endtrivlist}

\def\bs{{\bigskip}}

\begin{document}

\centerline {\bf Some Bounds for the Number of Components }
\smallskip \centerline {\bf of Real Zero Sets of Sparse Polynomials}

\bigskip

\centerline {Daniel Perrucci\footnote{Partially supported by the
Argentinian grant UBACyT 01-X198 and by CONICET}}

\bigskip

\centerline{Departamento de Matem\'atica, Facultad de Ciencias
Exactas y Naturales}

\centerline{Universidad de Buenos Aires}

\centerline {Ciudad Universitaria, Pab. I, (1428), Buenos Aires,
Argentina}

\centerline {e-mail address: perrucci@dm.uba.ar}

\begin{abstract}

We prove that the zero set of a $4$-nomial in $n$
variables in the positive orthant 
has at most three connected components. This bound, which does not
depend on the degree of the polynomial, not only improves the best
previously known bound (which was $10$) but is optimal as well.
In the general case, we prove that the number of connected
components of the zero set
of an $m$-nomial in $n$ variables in the positive orthant 
is lower than or equal to $(n+1)^{m-1}2^{1 + (m - 1)(m - 2)/2}$,
improving slightly the known bounds. 
Finally, we show that for generic exponents, the number of
non-compact connected components of the zero set of a $5$-nomial
in three variables in the positive octant is at most $12$. This
strongly improves the best previously known bound, which was
$10384$. All the bounds obtained in this paper continue to hold
for real exponents.

\end{abstract}

\section{Introduction.}

Descartes' Rule of Signs provides a bound for the number of
positive roots of a given real univariate polynomial which
depends on the number of sign changes among its coefficients but
not on its degree. One of its consequences is that the number of
positive roots of a polynomial with $m$ monomials is bounded
above by $m-1$.

Many attempts have been made to generalize Descartes' Rule of
Signs (or its corollaries) to a larger class of functions. Even
though this task has not yet been completed, important advances
have been made (\cite{GZ}, \cite{Kho1}, \cite{TesisKho}, \cite{Vor}, etc). 

Let us introduce the notation and terminology we
will use throughout this paper.
As usual, $\N$ will denote the set of positive integers. Let $n
\in \N$. Given $x \in \R_+^n = \{(x_1, \dots, x_n) \in \R^n \; |
\; x_k
> 0 \hbox { for } 1 \le k \le n\}$ and $a:= (a_1, \dots, a_n) \in \R^n$, $x^a$ will
denote $x_1^{a_1} \cdots x_n^{a_n}$.

\begin{defn}
Let $m \in \N$. An {\em $m$-nomial} in $n$ variables is a
function $f: \R_+^n  \to \R $ defined as
$$ f(x) = \sum_{i = 1}^m c_ix^{a_i},$$
where $c_i \in \R, c_i \ne 0$ and $a_i = (a_{i1}, a_{i2}, \dots,
a_{in}) \in \R^n$ for $i = 1, \dots, m$.
\end{defn}

An interesting fact is that \DR \ continues to hold if one counts multiplicities and
also if one allows real
exponents (the adaptation of the proof given in \cite[Proposition
1.1.10]{BR}, for instance, is straightforward).

\begin{defn}
\label{defK} Let $n, m \in \N$. Let us consider the functions $F:
\R_+^n \to \R^n$ of the form $F = (f_1, \dots, f_n)$ with $f_i$ an
$m_i$-nomial, 
such that the total number of distinct exponent vectors in $f_1,
\dots, f_n$ is less than or equal to $m$. We then define $K(n,m)$
to be the maximum number of isolated zeros (in $\R_+^n$) an
$F$ of this type may have.  
Similarly, we define $K'(n,m)$ to be the maximum number of
non-degenerate zeros (in $\R_+^n$) 
an $F$ of this type may have.
\end{defn}

A proof of the finiteness of $K(n,m)$ can be found in many
sources, for instance \cite[Corollary 4.3.8]{BR}. The finiteness
of $K'(n,m)$ is a consequence of the fact that $K'(n,m)$ is always
less than or equal to $K(n,m)$. A bound for $K'(n,m)$ is provided
by Khovanski's theorem, which is the most important result in the
theory of fewnomials:

\begin {thm}
\label{Khovanski} Following the notations above,
$$K'(n,m) \le (n+1)^{m-1}2^{(m-1)(m-2)/2}.$$
\end{thm}

For a proof of Khovanski's theorem, see \cite[Chapter 4]{BR},
\cite{Kho1} or \cite{TesisKho}. 
Nevertheless, the statement mentioned above is not exactly equal
to any of those in the references. To prove Theorem
\ref{Khovanski} divide every equation in the system $F(x) = 0$ by
$x^{a}$, where $x^a$ ($a \in \R^n$) is one of the monomials of
the system, to make the number of monomials drop and then use
\cite[Theorem 4.1.1]{BR} or \cite[Section 3.12, Corollary
6]{TesisKho}. Another fact to be highlighted is that here we are
allowing fewnomials with real exponents instead of integer
exponents as in \cite{BR}. Nevertheless, the proof in the last
reference does not make use of this fact.

Another way to generalize Descartes' Rule of Signs is to increase
just the number of variables. 
In this case, the problem is to find a bound for the number of
connected components of the zero set of a single polynomial, which
is expected to be a hypersurface. This paper is devoted to the
study of this problem, both in particular cases and in the general
one. The results presented here are inspired in a paper by Li,
Rojas and Wang (see \cite{LiRojWang}).

\begin{defn}
\label{TCN} Given a subset $X$ of $\R_+^n$, we will denote by
$\mathrm{Tot}(X)$, $\mathrm{Comp}(X)$ and $\mathrm{Non}(X)$ the
number of connected components, compact connected components and
non-compact connected components of $X$ respectively.

Given $n, m\in \N$, $P(n, m)$, $P_{comp}(n, m)$ and $P_{non}(n,
m)$ are defined in the following way. First, we define the set
$$\Omega(n,m) := \{f: \R_+^n \to \R \; |
\; f \hbox{ is a $k$-nomial with } 1 \le k \le m\}.$$ We then
define
$$P(n,m) := \max\{\, \mathrm{Tot}(f^{-1}(0))  \; | \; f \in \Omega(n,m)\},$$
$$P_{comp}(n,m) := \max\{\, \mathrm{Comp}(f^{-1}(0)) \; | \; f \in \Omega(n,m)\},$$
$$P_{non}(n,m) := \max\{\, \mathrm{Non}(f^{-1}(0)) \; | \; f \in \Omega(n,m)\}.$$

\end{defn}
It is clear from the definitions that, for all $n, m \in \N$,
$$P_{comp}(n,m) \le  P(n,m), \quad P_{non}(n,m) \le  P(n,m),$$
$$P(n,m) \le P_{comp}(n,m) + P_{non}(n,m)$$ and that $P$,
$P_{comp}$ and $P_{non}$ are increasing functions of their second
parameter. For fixed $n, m \in \N$, the finiteness of $P(n,m)$
(and thus that of $P_{comp}(n,m)$ and $P_{non}(n,m)$) is a
consequence of the fact that it is bounded from above by
$n(n+1)^m$ $ 2^{n-1}2^{m(m-1)/2}$ (see \cite{LiRojWang}, Corollary
2). Strongly based on this paper, we will derive a slightly better
bound:
\begin{thm} Using the previous notation,
\label{CuentaFinal} $P(n,m) \le (n+1)^{m-1}2^{1 + (m -
1)(m-2)/2}.$
\end{thm}

Our approach is different from that in \cite{LiRojWang} in the
way we bound the number of non-compact connected components.
We state our result in the following theorem, which will also be
useful in the last section, when dealing with $5$-nomials:

\begin{thm}
\label{Anexo} Let us consider $m,n \ge 2$. If $Z := f^{-1}(0)
\subset \R_+^n$ with $f$ an $m$-nomial in $n$ variables such that
the dimension of the Newton polytope (see Definition
\ref{Newton}) of $f$ is $n$, then
\begin{itemize}
\item $\mathrm{Non}(Z) \le 2n\,P(n-1,m-1),$ \item $\mathrm{Tot}(Z)
\le \sum_{i = 0}^{n-1} 2^i \frac{n!}{(n-i)!}P_{comp}(n-i,m-i).$
\end{itemize}
\end{thm}

Let us remark that, due to the fact that $\R_+^n$ is not a closed
set, a bounded connected component of the zero set of an
$m$-nomial may be non-compact. This is the case, for example, when
$f$ is the $3$-nomial in two variables defined by $f(x_1, x_2) =
x_1^2 + x_2^2 - 1.$

The next proposition shows that, for a fixed number of monomials, a big number of
variables will not increase the number of connected components:

\begin{prop} \label{nocrece} (see \cite[Theorem 2]{LiRojWang})
Given $m \in \N$, for all $n \in \N,$
$$P(n, m) \le  \begin{cases}
m-1 \hbox{\quad if } m \le  2, \cr P(m-2,m) \hbox{\quad if } m \ge
3. \cr
\end{cases}$$
\end{prop}

The reference given for the proposition above makes the additional assumption that $m \le
n+1$. Nevertheless, the proof there does not make use of this fact. As 
we really need to eliminate this extra assumption, 
we will give a brief proof of this proposition in the next section.  

One of the goals of this paper is to find a sharp bound for $P(n,
4)$ and Proposition \ref{nocrece} shows that it is enough to find such a
bound for $P(2,4)$. Our result is stated in the following theorem:

\begin{thm}
\label{Resfinal} Under the previous notation, we have:
\begin{enumerate}
\item $P_{comp}(2,4) = 1.$ \item $P_{non}(2,4) = 3.$ \item $P(2,4)
= 3$ (and thus $P(n,4) = 3$). \item If $f$ is a $4$-nomial in two
variables and $\dim \mathrm{Newt}(f) = 2$, then
$\mathrm{Tot}(f^{-1}(0)) \le 2.$
\end{enumerate}
\end{thm}

This theorem improves the best previously known bound for $P(n,
4)$,
which was $10$ (\cite[Theorems 2 and 3, and Example 2]{LiRojWang}). 
We will state the results used to prove this last bound and
sketch a brief proof of it in the next section.  Let us remark
that in \cite[Theorem 3]{LiRojWang} the equality of the second
item is proved in the  smooth  case.

The techniques we use to prove the previous theorems will also allow
us to prove the following theorem concerning $5$-nomials.

\begin{thm}
\label{five} Let $f$ be a $5$-nomial in three variables such that
$\dim \mathrm{Newt}(f) = 3$. Let $Z := f^{-1}(0) \subset \R^3_+.$
Then, $\mathrm{Non}(Z) \le 12$.
\end{thm}

This theorem significantly improves the best previously known
bound of $10384$ (the proof of this bound will be sketched briefly
in the next section too).

This paper is organized as follows: Section 2 details some
preliminaries. Section 3 concerns $4$-nomials and contains the
proof of Theorem \ref{Resfinal}. In Section $4$, we deal with the
general case of $m$-nomials in $n$ variables and we prove Theorems
\ref{CuentaFinal} and \ref{Anexo}. Finally, in Section 5, we
prove Theorem \ref{five}.

\section{Preliminaries}

\subsection{Previously known bounds for some particular cases}

The following result provides us with a bound for the
number of non-degenerate roots in the positive quadrant for a
fewnomial system having at most four different monomials.
\begin{lem} (See \cite[Section 2, Proposition 1]{LiRojWang}) Following the notation of Definition \ref{defK},
$K'(2,4) \le 5$.
\end{lem}

The next theorem enables us to get a bound for the number of
connected components in the positive orthant of the zero set of a
single fewnomial.
\begin{thm} \label{TeoLiRojWang} (see \cite[Theorem 2]{LiRojWang}) Following
the notation of Definition \ref{TCN}, we have:
\begin {itemize}
\item $P_{comp}(n,m) \le 2\lfloor K'(n,m)/2 \rfloor
 \le K'(n,m)$, 
\item $P_{non}(n,m) \le
2P(n-1,m)$.
\end {itemize}

\end{thm}

With these results, we can easily prove that $P(n,4) \le 10$ in
the following way:
$$P(n, 4) \le P(2, 4) \le P_{comp}(2, 4) + P_{non}(2, 4) \le
2\lfloor K'(2,4)/2\rfloor + 2P(1,4) \le 10,$$ the last inequality being true
because of \DR. We will improve this bound in Section 3.

In the same way,
$$P_{non}(3, 5) \le 2P(2, 5) \le 2P_{comp}(2, 5) + 2P_{non}(2, 5)
\le$$ $$\le 4\lfloor K'(2,5)/2\rfloor + 4P(1,5) \le 10384.$$ We will improve
this bound for the generic case in Section 5.

\subsection{Monomial changes of variables and Newton polytopes.}

Let us start this section with some notation and definitions.

\begin{notn}
Given a non-singular matrix $B \in \R^{n \times n}$, $B =
(b_{ij})_{1 \le i, j \le n},$ we will denote by $B_1$, $\dots$,
$B_n$ the columns of $B$. We will call the {\em monomial change
of variables associated to $B$} the function
$$h_B: \R_+^n \longrightarrow \R_+^n, \quad h_B(x) =
(x^{B_1}, \dots, x^{B_n}).$$
\end{notn}

The following formulae hold for all $x\in \R_+^n,$ $a \in \R^n$
and non-singular matrices $B, C \in \R^{n \times n}$:
\begin{itemize}
\item $h_B(x)^a = x^{Ba}.$
\item $ h_B \circ h_C =  h_{CB}.$
\end{itemize}

Recall the Newton polytope of a polynomial $f$, denoted by
$\mathrm{Newt}(f)$, which is a convenient combinatorial encoding
of the monomial term structure of a polynomial.

\begin{defn}
\label{Newton} Given an $m$-nomial $f$ in $n$ variables, $f(x)
:=\sum_{i=1}^m c_ix^{a_i}$, $\mathrm{Newt}(f)$ denotes the
smallest convex set containing the set of exponent vectors
$\{a_1, \dots ,a_m\}$. The dimension of $\mathrm{Newt}(f)$, $\dim
\mathrm{Newt}(f)$, is defined as the dimension of the smallest
translated linear subspace containing $\mathrm{Newt}(f)$.

\end{defn}

Therefore, for any $n$-variate $m$-nomial $f$, $\dim
\mathrm{Newt}(f) \le \min\{m-1, n\}.$

Given an $m$-nomial $f(x) = \sum_{i=1}^m c_ix^{a_i}$ and a
non-singular matrix $B\in \R^{n \times n}$, we have that
$$f \circ h_B(x) = \sum_{i=1}^m c_ih_B(x)^{a_i} = \sum_{i=1}^m c_ix^{Ba_i},$$
and thus
\begin{enumerate}
\item $f \circ h_B$ is also an $m$-nomial. \item $\mathrm{Newt}(f
\circ h_B) = \{B \, v \in \R^n \; | \; v \in \mathrm{Newt}(f)\}$,
and then, as $B$ is non-singular, $\dim \mathrm{Newt}(f) = \dim
\mathrm{Newt}(f \circ h_B)$.
\item As $h_B$ is an analytic automorphism of the positive orthant, then the
zero sets of $f$ and $f \circ h_B$ have the same number of compact and non-compact
connected components and critical points. 
\end{enumerate}

\begin{rem}
\label{sepuededividir} Given an $m$-nomial $f$ in $n$ variables,
$c \in \R, c \ne 0$ and $b \in \R^n$, the function $c^{-1}x^{-b}
f$ is an $m$-nomial whose Newton polytope is a translation of
$\mathrm{Newt}(f)$. Then $\dim \mathrm{Newt}(c^{-1}x^{-b}f) = \dim
\mathrm{Newt}(f)$. On the other hand, the zero set of
$c^{-1}x^{-b}f$ (included in $\R_+^n$ by definition) is equal to
the zero set of $f$ (also included in $\R_+^n$). In particular, by
choosing $c$ as one of the coefficients of $f$, we will get an
$m$-nomial with a coefficient equal to $1$. Moreover, by choosing
$b$ as one of the exponents of $f$, we will get an $m$-nomial with
a non-zero constant term. So, these particularities can be assumed
without loss of generality and not modifying the zero set of the
$m$-nomial, or the dimension of its Newton polytope. It can also
be proved that $p \in f^{-1}(0) \subset  \R^n_+$ is a critical
point of $f$ if and only if it is a critical point of
$c^{-1}x^{-b}f$.
\end{rem}

\begin{prop}
\label{dimcorr} Let $f$ be an $m$-nomial in $n$ variables, $Z :=
f^{-1}(0) \subset \R_+^n$ and $d := \dim \mathrm{Newt}(f).$ Then:
\begin{enumerate}
\item If  $d \le n-1$, then $\mathrm{Comp}(Z) = 0$ and
$\mathrm{Non}(Z) \le P(d,m).$ \item If $d = m - 1$, then
$\mathrm{Comp}(Z) = 0$ and $\mathrm{Non}(Z) \le 1.$
\end{enumerate}
\end{prop}

\begin{proof}{Proof:}
The proof can be done exactly as the proof of \cite[Theorem 2, Assertion 1]{LiRojWang}.
For instance, to prove the first assertion, suppose $f(x) = c_1
+ \sum_{i = 2}^{m}c_ix^{a_i}$. Let us consider a non-singular matrix $B \in
\R^{n \times n}$, such that the first $d$ columns of $B^{-1}$ are a basis of $\langle
a_2, \dots, a_m \rangle$. As each of the vectors
$Ba_i, i = 2, \dots, m,$ has its $n-d$ last coordinates
equal to zero, the $m$-nomial $f \circ h_{B}$
actually involves only $d$ variables
and its zero set may be described as
$Z' \times \R_+^{n-d}$, where $Z'$ is the zero set of an
$m$-nomial in $d$ variables. Thus, $\mathrm{Comp}(Z) = 0$ and
$\mathrm{Non}(Z) \le P(d, m)$.
\end{proof}

We can now give a proof of Proposition \ref{nocrece}.

\begin{proof}{Proof:}
Let $f$ be an $m$-nomial in $n$ variables, $Z := f^{-1}(0)$
$\subset \R_+^n$ and $d := \dim \mathrm{Newt}(f)$. For $m \le 2$, the proof is easy. 
If $m \ge 3$, as $d$ is always less than or equal to $m-1$, then we just need to consider the 
following cases:
\begin{itemize}
\item If $1 \le n \le m-2$, then $\mathrm{Tot}(Z) \le P(m-2,m)$,
because an $m$-nomial in $n$ variables can be considered as an
$m$-nomial in $m-2$ variables with the particularity that the last
$m - 2 - n$ variables are not actually involved in its formula.
\item If $m - 1 \le n$ and $d \le m - 2$, then $d \le n-1$. By
Proposition \ref{dimcorr}, $\mathrm{Tot}(Z) \le 0 + P(d,m) \le
P(m-2, m)$. \item If $m - 1 \le n$ and $d = m - 1$, again by
Proposition \ref{dimcorr}, $\mathrm{Tot}(Z) \le 0 + 1 \le P(m-2,
m)$.
\end{itemize}
\end{proof}

Finally, let us recall two classical results from topology that
will be quite useful in the next section.

\begin{thm}{(Connected curve classification.)}
Let $\Gamma$ be a differentiable manifold of dimension $1$. Then,
$\Gamma$ is  diffeomorphic either to $S^1$ or to  $\R$ depending
on whether $\Gamma$ is compact or not.
\end{thm}

The proof of this theorem can be found in  \cite{Milnor2}.

%% --------------------- Lema de Jordan y adaptacion ----------------------------

We will also use the next adaptation of Jordan's Lemma to the
positive quadrant, which can be easily proved from its original
statement (see, for example, \cite{Solom}) upon an application of
the exponential function.

\begin{lem}{(Adaptation of Jordan's Lemma).}
Let $\Gamma$ be a curve in $\R^2_+$ homeomorphic to $S^1$. Then,
$\R^2_+ \setminus \Gamma$ has two connected components, which we
will call $\mathrm{Int}(\Gamma)$ and $\mathrm{Ext}(\Gamma)$, such
that they are both open sets, $\mathrm{Int}(\Gamma)$ is bounded,
$\overline{\mathrm{Int}(\Gamma)} = \mathrm{Int}(\Gamma) \cup
\Gamma$ is compact and $\mathrm{Ext}(\Gamma)$ is unbounded.
\end{lem}

\section{On $4$-nomials in two variables. }

Most of the results we will obtain in this section come from the
study of the restriction of $4$-nomials in two variables to curves
of the type $\{ x \in \R_+^2 \; | \; x^a = J\}$ with $a \in \R^2$
and $J \in \R_+$. Let us introduce the notation we will use.

\begin{notn}
Let $f: \R_+^2 \to \R$ be an $m$-nomial in two variables, $p =
(p_1, p_2) \in \R_+^2$ and $u = (u_1, u_2)\in \R^2, u \ne 0$. By
$h_{(p, u)}$ we will denote the following parametrization of $\{
x \in \R_+^2 \; | \; x^u = p^u\}$:
$$h_{(p,u)}: \R_+ \to \R_+^2,$$
$$h_{(p, u)}(t) = (h^{(1)}_{(p, u)}(t), h^{(2)}_{(p, u)}(t)) = \left 
\{\begin{array}{ll}
(t, (p^u)^{1/u_2}t^{-u_1/u_2}) & \hbox{ if } u_2 \ne 0,\cr
(p_1, t) & \hbox{ if } u_2 = 0.
\end{array}
\right .$$
By $f_{(p,u)}$ we will denote the following function:
$$f_{(p, u)}: \R_+ \to \R, \quad f_{(p,u)} = f \circ h_{(p,u)}.$$
\end{notn}

\begin{rem}
\label{PtoCritSolit} \
\begin{itemize}
\item If $u_2 \ne 0$, $h_{(p,u)}(p_1) =
p$ and if $u_2 = 0$, then $h_{(p,u)}(p_2) = p$. \item $f_{(p, u)}$
is an $m'$-nomial in $1$ variable, with $m' \le m$. The
exponents of $f_{(p, u)}$ are proportional to 
the projections of the exponent vectors of $f$ on $\langle u
\rangle^\perp$. For instance, if $u_2 \ne 0$ and $a = (a_1, a_2)$
is an exponent of $f$, then $a_1 - u_1a_2/u_2 =  {\langle a, (u_2,
-u_1)\rangle}u_2^{-1}$ is an exponent vector of $f_{(p,u)}$. The
inequality $m'\le m$ is due to the fact that different exponent
vectors of $f$ may have the same projection on $\langle u
\rangle^\perp$, and so, some monomials in $f_{(p, u)}$ may be
re-grouped together and make the number of monomials decrease.
\item Suppose $p = (p_1, p_2)$ is a critical point of $f$
satisfying $f(p) = 0$ and $u = (u_1, u_2)\in \R^2$. If $u_2 \ne
0$, then $p_1$ is a degenerate zero of $f_{(p,u)}$, and if $u_2 =
0$, then $p_2$ is a degenerate zero of $f_{(p,u)}$. This is a
consequence of the chain rule.

\end{itemize}
\end{rem}

Notice that for $p \in \R_+^2$ and $u \in \R^2$, $u \ne 0$, the image of
$h_{(p,u)}$ is an unbounded curve containing $p$. The following
lemma will give us some information about the intersection between
this curve and a compact connected component of the zero set of an
$m$-nomial.

\begin{lem}
\label{CCSolit} Let $f$ be an $m$-nomial in two variables and let
$Z := f^{-1}(0) \subset \R^2_+$. Let $\Gamma$ be a compact
connected component of $Z$ containing only regular points of $f$
(so $\Gamma$ is a differentiable submanifold of $\R_+^2$
diffeomorphic to $S^1$). Let $p = (p_1, p_2) \in
\mathrm{Int}(\Gamma)$ and $u = (u_1, u_2)\in \R^2, u \ne 0$. Then,
if $u_2 \ne 0$, $f_{(p,u)}$ has a zero $s_1 \in (0, p_1)$ and a
zero $s_2 \in (p_1, + \infty)$ such
 that  $h_{(p,u)}(s_i) \in \Gamma$ $(i=1,2)$. If $u_2 = 0$,
$f_{(p,u)}$ has a zero $s_1 \in (0, p_2)$ and a zero $s_2 \in
(p_2, + \infty)$ such that  $h_{(p,u)}(s_i) \in \Gamma$ $(i=1,2)$.
\end{lem}

\begin{proof}{Proof:}

\bs

\centerline{\epsfysize 6cm\epsffile{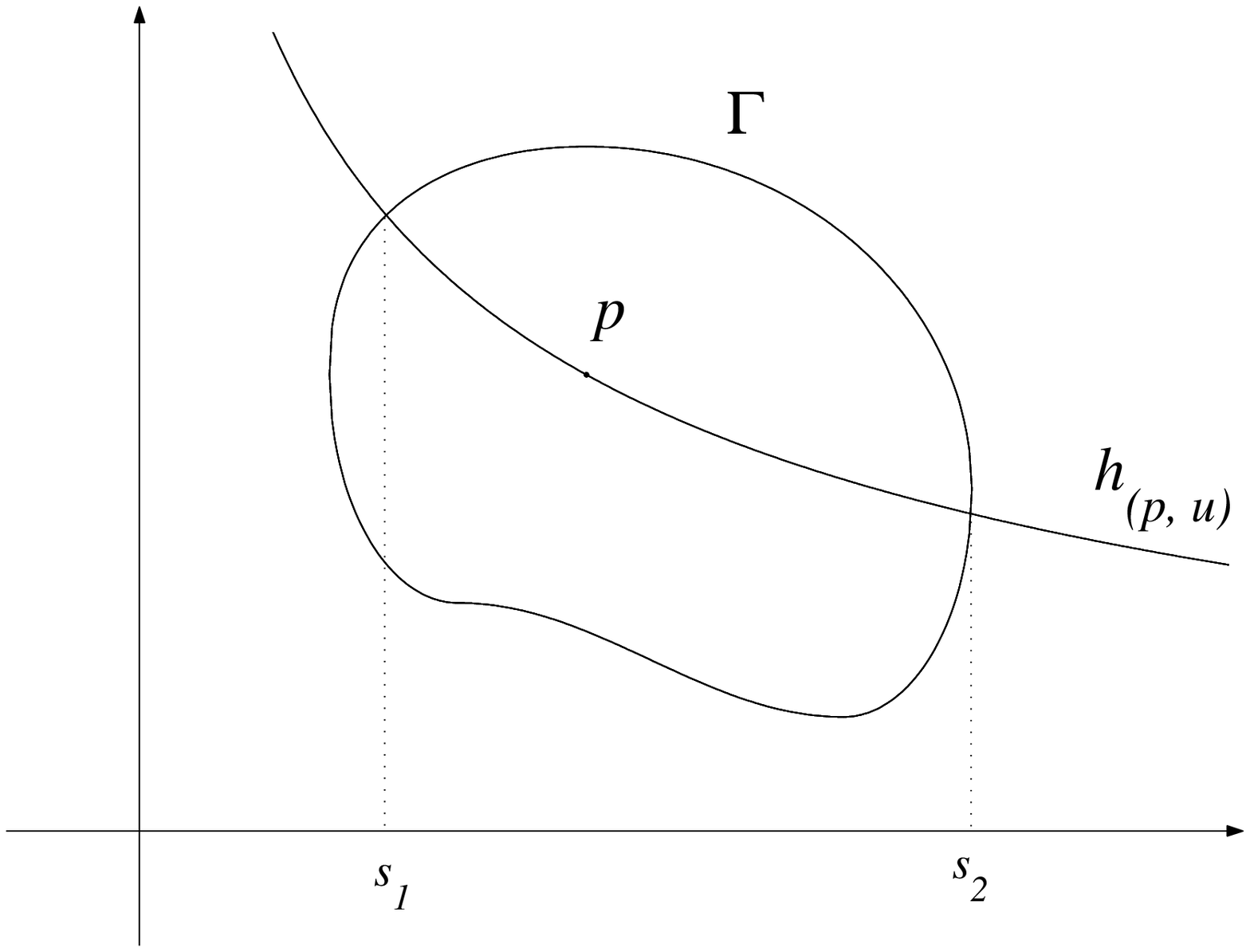}}

\bs

Let us suppose $u_2 \ne 0$. As $\Gamma$ is a compact set and $p$ lies in
$\mathrm{Int}(\Gamma)$, there exist $x \in (0, p_1)$ and $y \in (p_1, \infty)$
such that both $h_{(p,u)}(x)$ and $h_{(p,u)}(y)$ lie in $\mathrm{Ext}(\Gamma)$ and the lemma follows.
If $u_2 = 0$, a similar argument works.
\end{proof}

Suppose now that $f$ is a $4$-nomial in two variables. As
explained before, by studying the restriction of $f$ to curves of
a certain type we will obtain some information about its
coefficients.

\begin{lem}
\label{SobranCosas} Let $f$ be a $4$-nomial in two variables and
$Z := f^{-1}(0) \subset \R^2_+$. Suppose that one of the following
two conditions is satisfied:
\begin{enumerate}
\item $Z$ has a critical point $p = (p_1, p_2)$ and $Z \setminus \{p\} \ne \emptyset.$
\item $Z$ has a compact connected component $\Gamma$ and
$Z \setminus \Gamma \ne \emptyset.$
\end{enumerate}
Then, two of the coefficients of $f$ are positive and the other
two are negative.
\end{lem}

\begin{proof}{Proof:}
Suppose $Z$ satisfies the first condition. Let $q = (q_1, q_2) \in
Z \setminus \{p\}$.

If $p_1 \ne q_1$, then $p_1/q_1 \ne 1$. Let
$$v_1 := \frac{\log(q_2/p_2)}{\log(p_1/q_1)}$$
Then $p_1^{v_1}p_2 = q_1^{v_1}q_2$. Let $v \in \R^2$, $v := (v_1,
1)$. As it was explained in Remark \ref{PtoCritSolit}, $p_1$ is a
zero of $f_{(p,v)}$ with multiplicity at least $2$. On the other
hand,
$$f_{(p,v)}(q_1) = f(q_1, p^vq_1^{-v_1}) = f(q_1, q^vq_1^{-v_1}) = f(q) = 0,$$
because $q \in Z.$ As $p_1 \ne q_1$, we know that $f_{(p,v)}$ has
at least three zeros (counting multiplicities) in $\R_+$. We know
that $f_{(p,v)}$ is an $m'$-nomial with $m' \le 4$. By Descartes'
Rule of Signs, we know that the number of sign changes in
$f_{(p,v)}$ is at least three; thus, $m' = 4$ and among the $4$
coefficients of $f_{(p,v)}$, there must be two positive and two
negative. On the other hand, if
$$f(x) = \sum_{i = 1}^4 c_ix^{a_i},$$
then
$$f_{(p,v)}(x_1) = \sum_{i = 1}^4 c_i \, (p^v)^{a_{i2}} x_1^{a_{i1}- a_{i2}v_1}.$$
As the signs of the coefficients of $f_{(p,v)}$ are defined by the
signs of the coefficients of $f$, then $f$ must have two positive
 and two negative coefficients.

If $p_1 = q_1$, as  $p \ne q$, we will have $p_2 \ne q_2$. In this
case, let us take $v := (1,0)$ and proceed as above. 

Let us suppose now that $Z$ satisfies the second condition, which
is having a compact connected component $\Gamma$, and $Z \ne
\Gamma$. If $Z$ has a critical point, then the first condition is
also satisfied. If it does not have a critical point, we consider
$\hat p := (\hat p_1, \hat p_2) \in \mathrm{Int}(\Gamma)$ and
$\hat q := (\hat q_1, \hat q_2) \in Z \setminus \Gamma$.

If $\hat p_1 \ne \hat q_1$, in the same way we did before, we can
find a vector $w \in \R^2$, $w = (w_1, 1)$ such that $\hat p^w =
\hat q^w$. Then, $f_{(\hat p,w)}$ has at least one zero $s_1$ in
the interval $(0, \hat p_1)$ such that $h_{(\hat p,w)}(s_1) \in
\Gamma$ and at least one zero $s_2$ in the interval $(\hat p_1, +
\infty)$ such that $h_{(\hat p,w)}(s_2) \in \Gamma$. On the other
hand,
$$f_{(\hat p,w)}(\hat q_1)
= f(    \hat q_1     ,       \hat p^w \hat q_1^{-w_1} ) = f( \hat
q_1     ,       \hat q^w \hat q_1^{-w_1}   ) = f(\hat q) = 0,$$
because $\hat q \in Z.$ Besides, due to the fact that $h_{(\hat
p,w)}(\hat q_1) = \hat q \in Z \setminus \Gamma,$ $\hat q_1 \ne
s_1$ and $\hat q_1 \ne s_2$. Then, we deduce that $f_{(\hat
p,w)}$ has at least three zeros in $\R_+$, and then $f_{(\hat
p,w)}$ is also a $4$-nomial with at least three sign changes. So,
$f_{(\hat p,w)}$ and $f$ have both two coefficients with each
sign.

If $\hat p_1 = \hat q_1$, then  $\hat p_2 \ne \hat q_2$, and the
same argument works.
\end{proof}

Due to the lemma above, we will focus our attention for a moment
on $4$-nomials with two coefficients of each sign. We will start
relating some properties of the zero set of a $4$-nomial in two
variables of this form with its Newton polytope.

\begin{lem}
\label{CuadriOpu} Let $f$ be a $4$-nomial in two variables with
two positive and two negative coefficients, such that $\dim
\mathrm{Newt}(f) = 2$. Let $Z := f^{-1}(0) \subset \R^2_+$, and
suppose one of the following two conditions is satisfied:
\begin{enumerate}
\item $Z$ has a critical point $p = (p_1, p_2)$.
\item $Z$ has a compact connected component $\Gamma$.
\end{enumerate}
Then $\mathrm{Newt}(f)$ is a quadrilateral without parallel sides
and 
coefficients corresponding to adjacent vertices have opposite
signs.
\end{lem}

\begin{proof}{Proof:}
Define $r = (r_1, r_2) \in \R^2_+$ as follows: if $Z$ satisfies
the first condition, then $r = p$ and if $Z$ satisfies the second
one but not the first one (so $\Gamma$ is diffeomorphic to
$S^1$), then $r$ is any point in $\mathrm{Int}(\Gamma).$ By Remark
\ref{PtoCritSolit} and Lemma \ref{CCSolit}, we know that for all
$v = (v_1, v_2) \in \R^2$, $v \ne 0$,  $f_{(r,v)}$ has at least
two zeros (counting multiplicities) in $\R_+$.

Since $\dim \mathrm{Newt}(f) = 2$, the exponent vectors do not lie
on a line. Suppose $f(x) = \sum_{i = 1}^4c_ix^{a_i}$. Then the
vertices of $\mathrm{Newt}(f)$ are among the vectors $a_1, a_2,
a_3$ and $a_4$ and $\mathrm{Newt}(f)$ might either be a triangle
or a quadrilateral. We will need to study  four cases separately.

\begin{itemize}
\item Suppose $\mathrm{Newt}(f)$ is a triangle whose vertices are
the vectors $a_1, a_2$ and $a_3$ and that the vector $a_4$ lies in
the interior of $\mathrm{Newt}(f)$. Assume $c_1$ and $c_2$ are
positive and $c_3$ and $c_4$ are negative (by multiplying $f$ by
$-1$ and reordering the monomials if necessary). Let $v := a_1 -
a_4 \ne 0$ and $L$ the line through $a_1$ and $a_4$.

\bs

\centerline{\epsfysize 4cm\epsffile{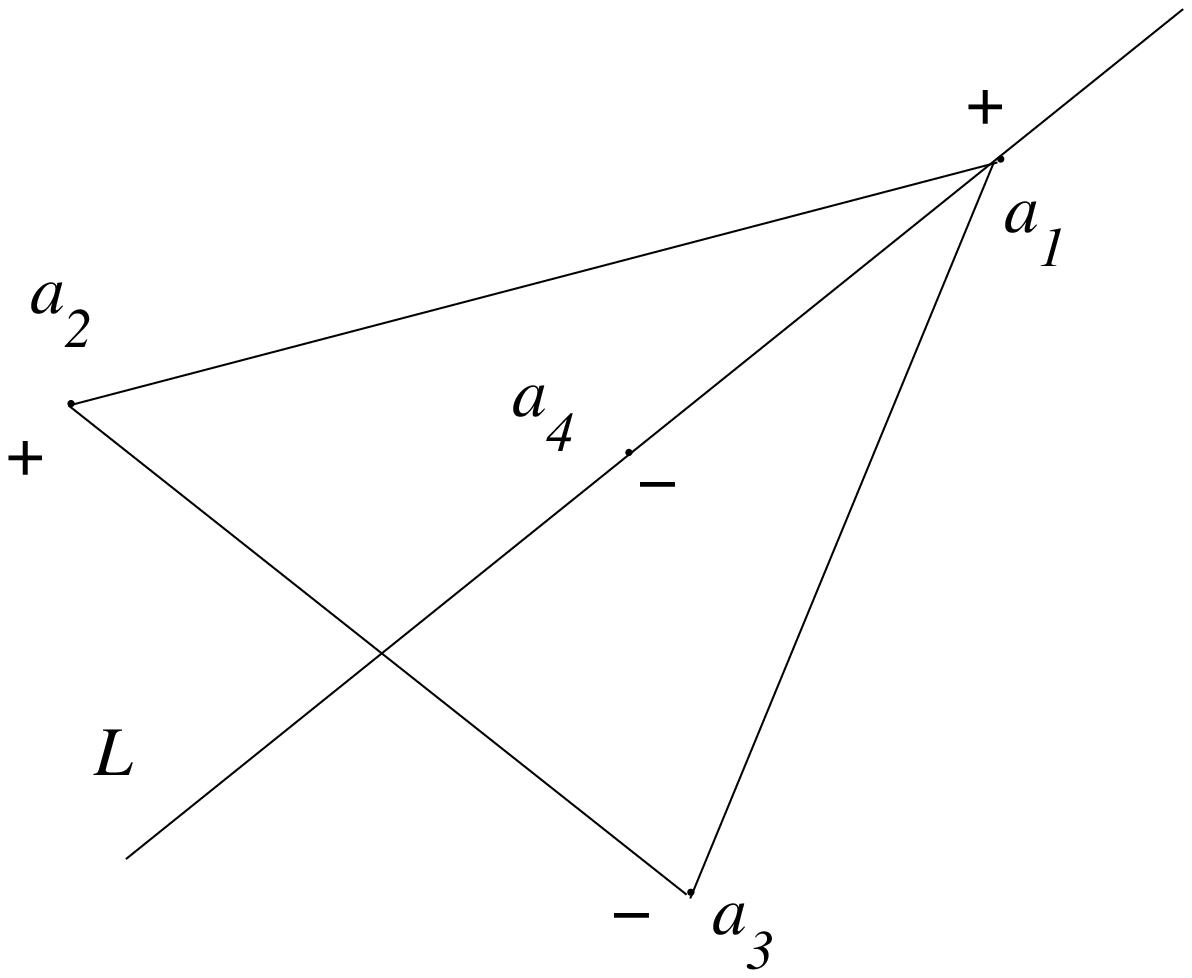}}

\bs

As $a_1$ and $a_4$ have the same projection on $\langle
v\rangle^\perp$ and $a_2$ and $a_3$ are on opposite sides of the
line $L$, we conclude that $f_{(r, v)}$ is a $3$-nomial of the
following type (if $v_2 \ne 0$):
$$f_{(r, v)}(x_1) = c_3(r^v)^{a_{32}/v_2} x_1^{a_{31}-
a_{32}v_1/v_2} + $$
$$ + (c_1(r^v)^{a_{12}/v_2} + c_4(r^v)^{a_{42}/v_2}) x_1^{a_{41}-a_{42}v_1/v_2} +
 c_2(r^v)^{a_{22}/v_2} x_1^{a_{21}- a_{22}v_1/v_2}.$$
Even though we do not know at this point if in the above formula
the terms are written in increasing or decreasing order, this
$3$-nomial has exactly one sign change, because the monomials of
higher and lower exponent have distinct coefficient signs (we are
supposing $c_2 > 0$ and $c_3 < 0$). By \DR, it cannot have $2$
zeros (counting multiplicities) as we know it does. Then, we have
a contradiction, and we conclude that the Newton polytope of $f$
cannot be a triangle having the remaining exponent vector in its
interior. If $v_2 = 0$, the same procedure works.

\item Let us suppose now that $\mathrm{Newt}(f)$ is a triangle whose
vertices are the exponent vectors $a_1, a_2$ and $a_3$; and that
the vector $a_4$ lies on one of the edges of $\mathrm{Newt}(f)$.
Without loss of generality, we suppose that $a_4$ lies on the
segment $a_1a_2$.

\bs

\centerline{\epsfysize 4cm\epsffile{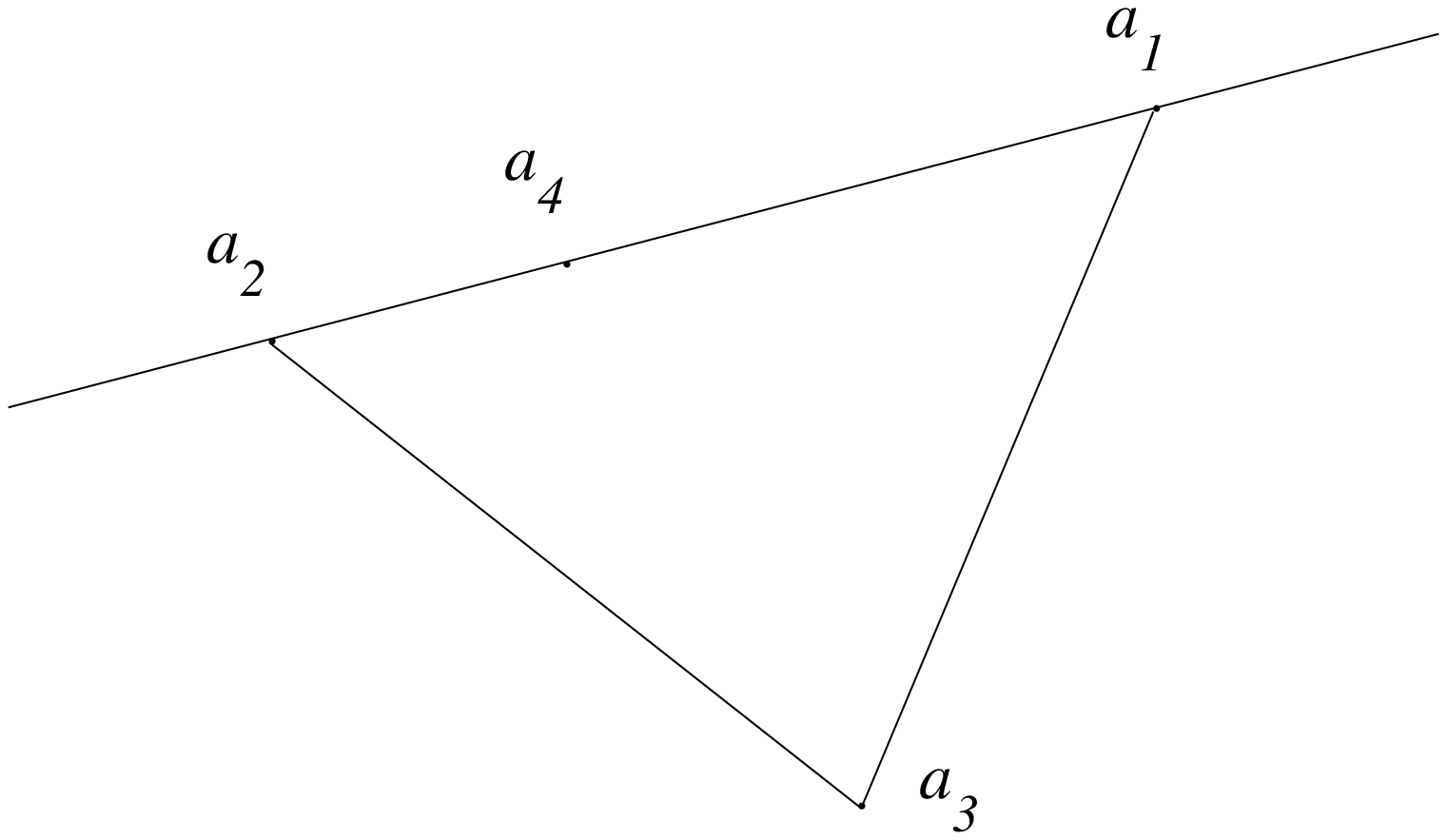}}

\bs

By taking again $v := a_1 - a_4$, we have that $a_1$, $a_2$ and
$a_4$ have the same projection on $\langle v\rangle^{\perp}$.
Thus, $f_{(r, v)}$ is a $2$-nomial (because its first, second and
fourth term can be re-grouped together in a single monomial) and,
by \DR, $f_{(r, v)}$ cannot have two zeros (counting
multiplicities) as we know it should. We have then  a
contradiction, which enables us to eliminate this case.

\item Suppose $\mathrm{Newt}(f)$ is a quadrilateral with a pair of
parallel opposite sides. Without loss of  generality, we suppose
that the segments $a_1a_2$ and $a_3a_4$ are parallel.

\bs

\centerline{\epsfysize 4cm\epsffile{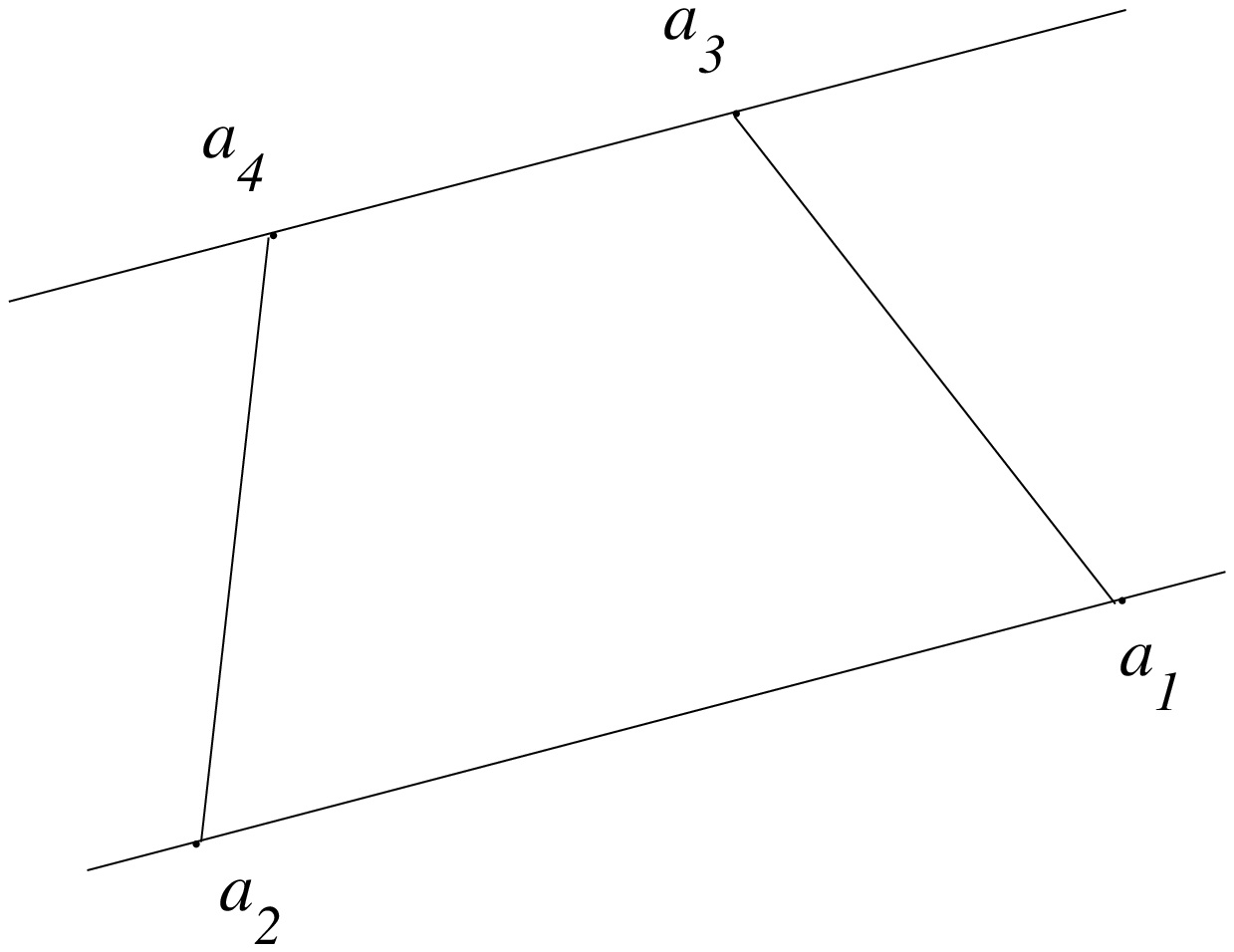}}

\bs

Let us take $v := a_1 - a_2$. As  $a_1$ and $a_2$ have the same
projection on $\langle v\rangle^\perp$, and $a_3$ and $a_4$ also
do so, we can re-group the monomials in $f_{(r, v)}$ and form a
$2$-nomial, which again is impossible.

\item Finally, suppose that $\mathrm{Newt}(f)$ is a quadrilateral and
that the coefficients of same sign correspond to adjacent
vertices. Without loss of generality, let us suppose that $a_1$
and $a_2$ are adjacent, $a_3$ and $a_4$ are adjacent too, $c_1$
and $c_2$ are positive and $c_3$ and $c_4$ are negative.

\bs

\centerline{\epsfysize 4cm\epsffile{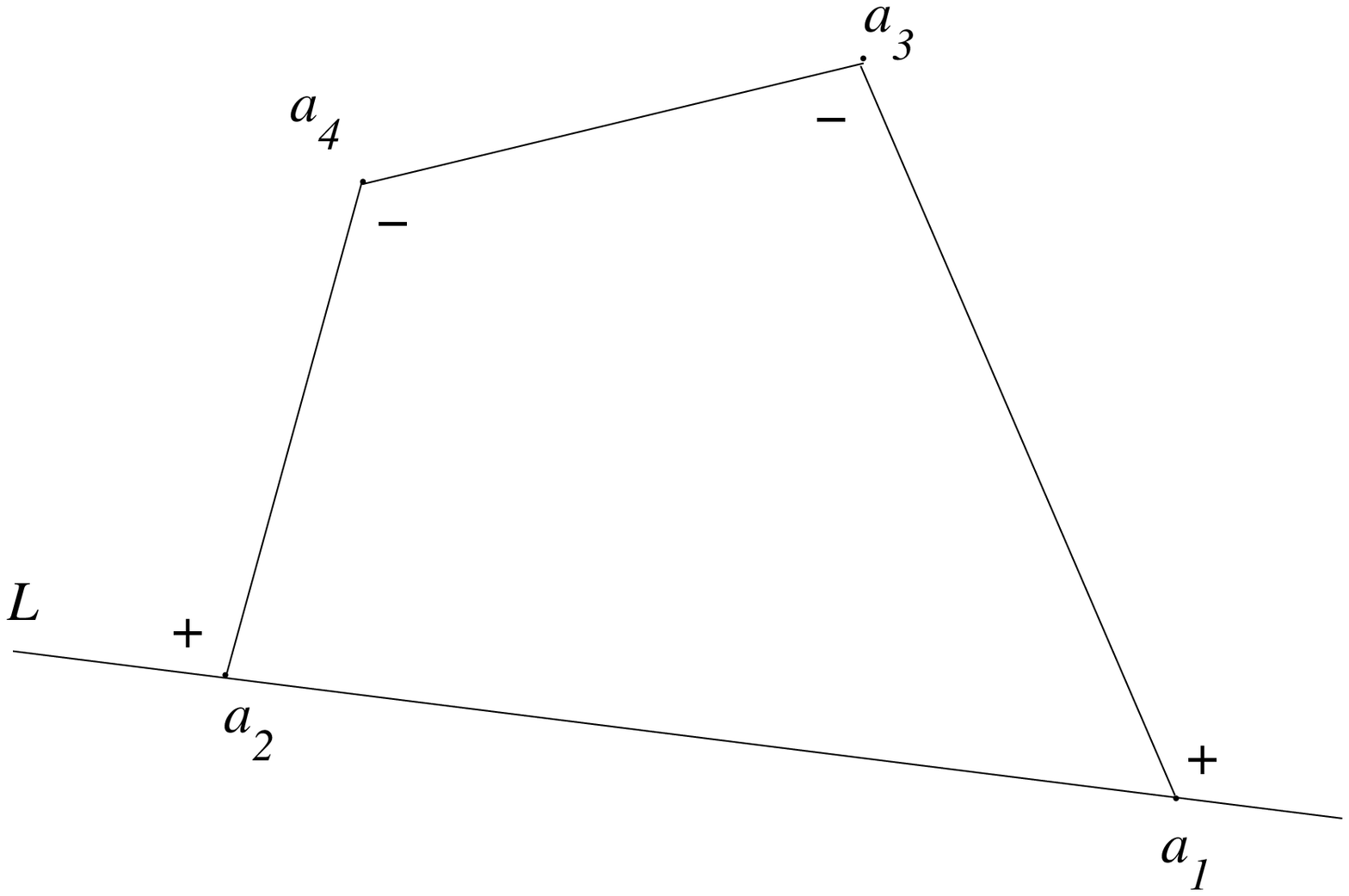}}

\bs

Let $v := a_1 - a_2$ and $L$ the line through $a_1$ and $a_2$. As
$a_1$ and $a_2$ have the same projection on $\langle
v\rangle^\perp$, then $f_{(r,v)}$ is a $3$-nomial. But, as the
two remaining exponent vectors (corresponding both to negative
coefficients) lie in the same side of $L$, $f_{(r,v)}$  has just
one sign change. For this reason, it cannot have two zeros,
and we get a contradiction. 
\end{itemize}
We conclude that the lemma follows.
\end{proof}

The next lemma shows the existence of a convenient change of
variables for certain bivariate $4$-nomials.

\begin{lem}
\label{Equiv} Let $f$ be a bivariate $4$-nomial having two
positive coefficients and two negative coefficients and such that
$\mathrm{Newt}(f)$ is a quadrilateral with no parallel opposite
sides and coefficients corresponding to adjacent vertices having
opposite signs. Then, there is an invertible change of variables
$h$ such that $f \circ h$ is
$$f \circ h(x_1, x_2) = 1 - x_1 - x_2 +
Ax_1^cx_2^d,$$ with $A > 0$, $c, d > 1$, and $h$ is the
composition of a monomial change of variables with a re-scaling of
the variables.
\end{lem}

\begin{proof}{Proof:}
Suppose that we enumerate the vertices of $\mathrm{Newt}(f)$ in
such a way that $a_1$ and $a_4$ are not adjacent. Because of
Remark \ref{sepuededividir} we can suppose $f$ is of the
following type:
$$f(x_1, x_2) = 1 + \sum_{i = 2}^4 c_ix^{a_i},$$
i.e., $a_1 = (0,0)$. As coefficients with the same sign correspond
to non-adjacent vertices of $\mathrm{Newt}(f)$, we know that $c_2,
c_3 <0$ and $c_4
> 0$. Consider the four triangles that can be formed with three
of the four vertices of $\mathrm{Newt}(f)$. Among these triangles,
there must be one having the minimal area. Suppose it is the
triangle $a_1a_2a_3$ (this can be enforced by rotating indices if
necessary). Because of the fact that $\mathrm{Newt}(f)$ does not
have parallel opposite sides, this area is strictly less than the
area of the triangles $a_1a_2a_4$ and $a_1a_3a_4$.

Let $B := \{a_2,a_3\}$ and $C$ be the matrix having the elements
of $B$ as columns. As $\mathrm{Newt}(f)$ is a quadrilateral, $a_1
= (0,0), a_2$ and $a_3$ do not lie on a line. Then $B$ is a basis
of $\R^2$ and $C$ is non-singular. As in \cite[Lemma
1]{LiRojWang}, let $h$ be the composition of $h_{C^{-1}}$ and the
linear re-scaling $(x_1, x_2) \mapsto \big( \frac{x_1}{|c_2|},
\frac{x_2}{|c_3|}\big)$. Let $(c, d) := C^{-1}a_4$. Then,
$$f \circ h(x_1, x_2) =
1 - x_1 - x_2 +c_4\frac{1}{|c_2|^c}\frac{1}{|c_3|^d} x_1^c
x_2^d.$$

Let $A$ be the last coefficient of the $4$-nomial above. Then $A
> 0$. On the other hand, the Newton polytope $\mathrm{Newt}(f \circ
h_{C^{-1}})$ must also be a quadrilateral with vertices $(0,0),
(1,0)$, $(0,1)$ and $(c,d)$. As $a_1$ and $a_4$ are opposite
vertices in $\mathrm{Newt}(f)$, then $(0,0)$ and $(c, d)$ must be
opposite vertices in $\mathrm{Newt}(f \circ h_{C^{-1}})$. Thus, $c
, d > 0$. As the area of triangle $a_1a_2a_3$ is smaller than
that of triangle $a_1a_2a_4$, the area of triangle $(0,0) (1,0)
(0,1)$ should be smaller than that of triangle $(0,0) (1,0)
(c,d)$, and thus $d > 1$. In an analogous way, we can prove that
$c> 1$.
\end{proof}

Let us recall that, as $h$ is a diffeomorphism of $\R_+^2$, the
zero sets of $f$ and $f \circ h$ have the same number of compact
and non-compact connected components and critical points.

The following lemma will let us deal with the case when the zero
set of the $4$-nomial has a critical point.

\begin{lem}
\label{PtoCrit} Let $f$ be a $4$-nomial in two variables such that
$\dim \mathrm{Newt}(f) = 2$ and let $Z := f^{-1}(0) \subset
\R^2_+$. Suppose that $p = (p_1, p_2) \in Z$ is a critical point
of $f$, and also that $Z \setminus \{p\} \ne \emptyset$. Then $Z$
is a connected non-compact set; that is to say, $\mathrm{Non}(Z)
= 1$ and $\mathrm{Comp}(Z) = 0$.
\end{lem}

\begin{proof}{Proof:}

The proof of this lemma will be done in four steps. In the first
one, we will make use of the lemmata
we proved before to make sure
that it is enough to restrict our attention to $4$-nomials of a
very specific form. In the second one, we will study the sign of
$f$ on some curves we will consider. In the third one, we will use
the information we obtained to characterize a non-compact
connected set $W$ where $f$ vanishes. In the last one, we will
prove that, in fact, $W = Z$.

{\it Step 1.}\ By Lemma \ref{SobranCosas}, we know that, among the
coefficients of $f$ there must be two positive and two negative
ones, and by Lemma \ref{CuadriOpu}, $\mathrm{Newt}(f)$ must be a
quadrilateral without parallel opposite sides, and with same sign
coefficients corresponding to opposite vertices. Then, by Lemma
\ref{Equiv}, we can suppose $f$ is of the following type:
$$f(x_1, x_2) = 1 - x_1 - x_2 + Ax_1^cx_2^d,$$
with $A > 0$ and $c, d > 1$.

{\it Step 2.}\ Let $v := (c,d-1)$. This vector has the nice
property that $(0,1)$ and $(c,d)$, which are exponent vectors in
$f$, have the same projection on $\langle v \rangle^{\perp}$.
Because of this fact, $f_{(p,v)}$ is a $3$-nomial. In fact,
$$f_{(p,v)}(x_1) =  \big(- (p^v)^{1/(d-1)} + A \, (p^v)^{d/(d-1)} \big) x_1^{-c/(d-1)} + 1 - x_1.$$
Notice that $-c/(d-1) < 0$ because $c, d > 1$. By Remark
\ref{PtoCritSolit}, $p_1$ is a zero of $f_{(p,v)}$ of multiplicity
greater than or equal to $2$. By \DR, the $3$-nomial $f_{(p,v)}$
must have at least two sign changes, $p_1$ is a zero of
multiplicity exactly $2$ and  $f_{(p,v)}$ does not have other
zeros. As the unique zero of $f_{(p,v)}$ has an even
multiplicity, and its leading exponent coefficient is negative,
we know that $f_{(p,v)}(x_1) \le 0$ for all $x_1 \in \R_+$ and
$f_{(p,v)}(x_1) < 0$ if $x_1 \ne p_1$.

As $(0,0)$ and $e_2:=(0,1)$ are both vector exponents in $f$, it
can be proved in an analogous way that $f_{(p,e_2)}(x_1) \ge 0$
for all $x_1 \in \R_+$ and $f_{(p,e_2)}(x_1) > 0$ if $x_1 \ne
p_1$. Moreover,
$$f_{(p,e_2)}(x_1) = f(x_1, (p^{e_2})^1x_1^0) =
(1 - p_2) -x_1 + Ap_2^dx_1^c,$$ and as it has a zero of
multiplicity equal to two, it has two sign changes and then we
deduce that $p_2 < 1$. In the same way we can also prove that
$p_1 < 1$.

It can easily be checked that for $x_1 \in (0, p_1),
h^{(2)}_{(p,e_2)}(x_1) < h^{(2)}_{(p,v)}(x_1)$, and for $x_1 \in
(p_1, + \infty), h^{(2)}_{(p,e_2)}(x_1) > h^{(2)}_{(p,v)}(x_1)$.
To illustrate the situation, in the following figure we have drawn
the curves $h_{(p, e_2)}$ and $h_{(p, v)}$ indicating the sign of $f$ on them:

\bs

\centerline{\epsfysize 7.4cm\epsffile{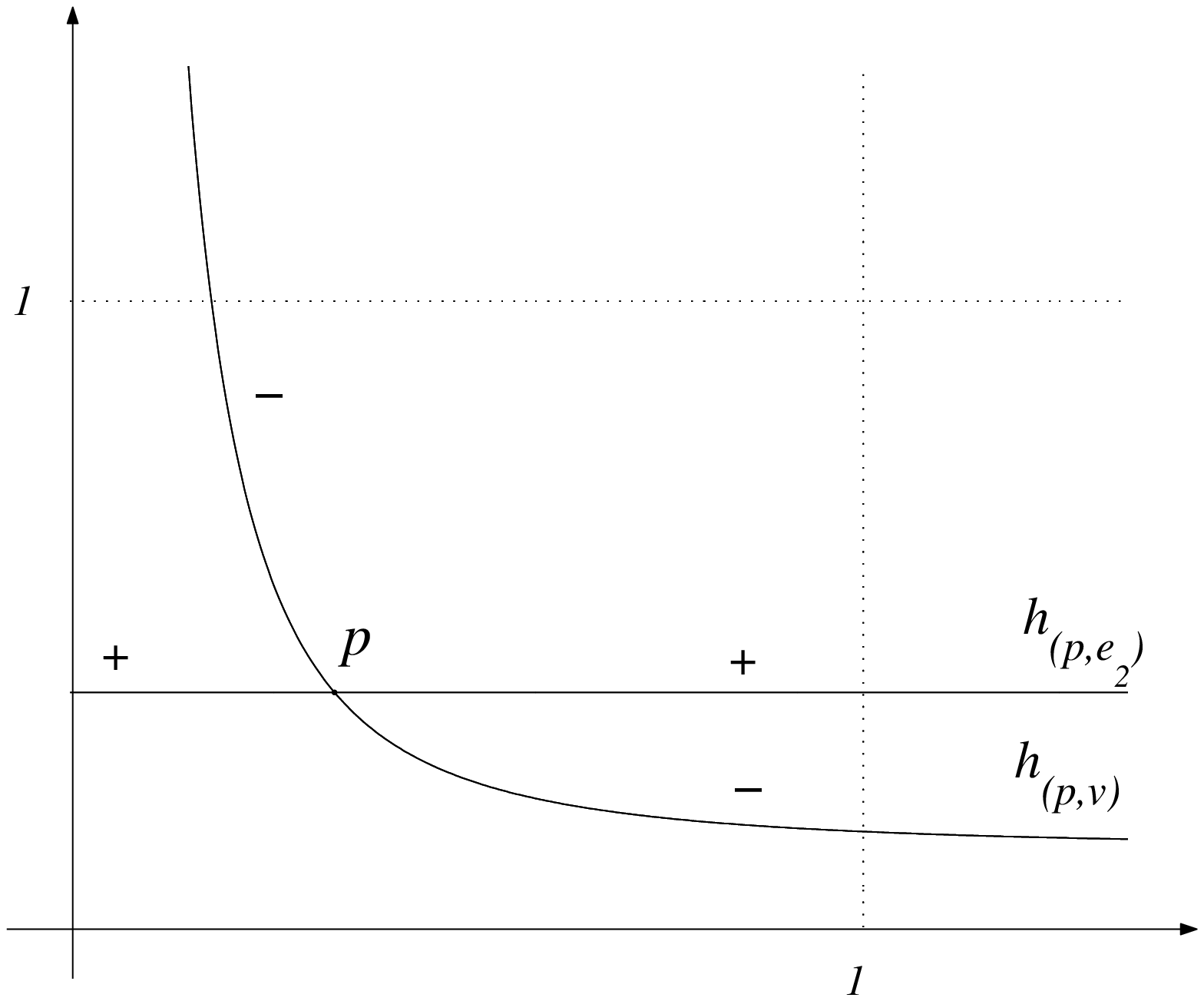}}

\bs

Finally, for a fixed $\alpha \in \R_+$, let us analyze the
function $f(\alpha,x_2)$ in the variable $x_2$:
$$f(\alpha,x_2) = (1 - \alpha) - x_2 + A\alpha^cx_2^d.$$
Let us notice that for every fixed $\alpha \in (0,1)$,
$\lim\limits_{x_2 \to 0+} f(\alpha, x_2) =
1 - \alpha > 0.$

{\it Step 3.}\ In order to study how many times the line $\{x_1 =
\alpha\}$ intersects $Z$ for a fixed $\alpha \in \R_+$, we will
continue studying the function $f(\alpha,x_2)$. As $A\alpha^c >
0$ and $d > 1$, if $\alpha < 1$, this function is a $3$-nomial
with two sign changes. Because of Descartes' Rule of Signs, it
will have either no zeros or two (counted with multiplicity) in
$\R_+$. If $\alpha = 1$, this function is a $2$-nomial with just
one sign change and, finally, if $\alpha
> 1$, it is a $3$-nomial with one sign change. In both cases,
it has exactly one zero in $\R_+.$

For a fixed $\alpha \in (0, p_1)$ the function (in the variable
$x_2$) $f(\alpha,x_2)$ 
must have an odd number of zeros (counted with multiplicity) in
the interval $(h^{(2)}_{(p,e_2)}(\alpha),$
$h^{(2)}_{(p,v)}(\alpha))$. As it has at most two zeros in
$\R_+$, then it has just one zero in that interval. Let us call it
$g(\alpha)$.

In an analogous way, for a fixed $\alpha \in (p_1, 1)$ the
function 
$f(\alpha,x_2)$ 
must have at least one zero in the interval
$(0,h^{(2)}_{(p,v)}(\alpha))$, which we will call $t(\alpha)$,
and another one in the interval $(h^{(2)}_{(p,v)}(\alpha),
h^{(2)}_{(p,e_2)}(\alpha))$, which we will call $g(\alpha)$. As
this function has at most two zeros in $\R_+$, then it has no
other zeros.

For a fixed $\alpha \in [1, +\infty)$, the function
$f(\alpha,x_2)$ must have an odd number of zeros (counted with
multiplicity) in the interval $(h^{(2)}_{(p,v)}(\alpha),
h^{(2)}_{(p,e_2)}(\alpha))$. As this function has at most one
zero in $\R_+$, then it has just one zero in that interval.
Again, let us call it $g(\alpha)$.

Finally, let us define $g(p_1) = p_2$, and let us prove that the
function $g:\R_+ \to \R_+$ we have just defined is continuous.

As $p$ is a critical point of $f$, we know that
$$\frac{\partial f}{\partial x_2}(p) = -1 + dAp_1^{c}p_2^{d-1} = 0,$$
and this implies that $p^v = 1/dA.$

Suppose there exists $x_1 \in \R_+, x_1 \ne p_1$ such that
$\frac{\partial f}{\partial x_2}(x_1, g(x_1)) = 0$; then
$$g(x_1) = (1/dA)^{1/(d-1)}x_1^{-c/(d-1)} = (p^v)^{1/(d-1)}x_1^{-c/(d-1)} =
h^{(2)}_{(p,v)}(x_1),$$ and this is impossible because of the
definition of $g$. Then, for all $x_1 \ne p_1$, $\frac{\partial
f}{\partial x_2}(x_1, g(x_1)) \ne 0$.

Let us fix $\alpha\ne p_1$ and see that $g$ is continuous in
$\alpha$. Suppose that $\alpha > p_1$ (if $\alpha < p_1$ the proof
can be done in the same way). We know that
$h^{(2)}_{(p,v)}(\alpha)= (p^v)^{1/d-1}\alpha^{-c/d-1} < g(\alpha)
< p_2 = h^{(2)}_{(p,e_2)}(\alpha)$ and $\frac{\partial f}{\partial x_2}(\alpha, g(\alpha)) \ne
0.$ Then, by the Implicit Function Theorem, there is a continuous
function, let us call it $s$, defined in an interval $(\alpha -
\varepsilon, \alpha + \varepsilon)$ with $\alpha - \varepsilon
> p_1$, such that
$ s(\alpha) = g(\alpha)$ and $f(x_1, s(x_1)) = 0$ for all $x_1$
in the interval of definition. Moreover, choosing a suitable value
of $\varepsilon$, we can suppose that, for all $x_1$ in $(\alpha
- \varepsilon, \alpha + \varepsilon)$, $s(x_1)$ lies in
$((p^v)^{1/(d-1)} x_1^{-c/(d-1)}, p_2)$. As $x_2 = g(x_1)$ is the
unique value in this interval such that $f(x_1, x_2) = 0$, we
have $g \equiv s$ in $(\alpha - \varepsilon, \alpha +
\varepsilon)$ and therefore $g$ is continuous in $\alpha$.

To prove that $g$ is continuous in $p_1$, notice that if $x_1 >
p_1$, then
$$ (p^v)^{1/(d-1)} x_1^{-c/(d-1)} < g(x_1) < p_2$$
and
$$\lim\limits_{x_1 \to p_1^+} (p^v)^{1/(d-1)} x_1^{-c/(d-1)} = (p^v)^{1/(d-1)} p_1^{-c/(d-1)} = p_2.$$
So, we have $\lim\limits_{x_1 \to p_1^+}g(x_1) = p_2$.
Analogously, we prove that $\lim\limits_{x_1 \to p_1^-}g(x_1) =
p_2$.

\smallskip

Now, let us consider $w \in \R^2, w := (c-1,d)$. In the same way
we proved the existence of the function $g$, we can prove that
there exists a function $k: \R_+ \to \R_+$ satisfying the
following properties:
\begin{itemize}
\item For all positive $x_2$, we have $f(k(x_2), x_2) = 0.$
\item $k$ is continuous.
\item If $x_2 < p_2$ then $k(x_2) \in (p_1, (p^w)^{1/(c-1)}x_2^{-d/(c-1)})$, and $x_1 = k(x_2)$
is the unique value in that interval such that $f(x_1, x_2) = 0.$
\item $k(p_2) = p_1$.
\item If $x_2 > p_2$ then $k(x_2) \in ((p^w)^{1/(c-1)}x_2^{-d/(c-1)}, p_1)$, and $x_1 = k(x_2)$
is the unique value in that interval such that $f(x_1, x_2) = 0.$
\end{itemize}

Let us define $W_1 = \{(x_1, g(x_1)) \; | \; x_1 \in \R_+\}
\subset \R_+^2$, $W_2 = \{(k(x_2), x_2) \; | \; x_2 \in \R_+\}
\subset \R_+^2$ and $W = W_1 \cup W_2$. As the functions $g$ and
$k$ are continuous, $W_1$ and $W_2$ are connected. As $g(p_1) =
p_2$ and $k(p_2) = p_1$, it follows that $p \in W_1 \cap W_2$, and
then $W$ is connected. Moreover, it is an unbounded set.

{\it Step 4.}\ Let us prove now that $W = Z$, and therefore, that
$\mathrm{Non}(Z) = 1$ and $\mathrm{Comp}(Z) = 0$.

Due to the fact that, for all $x_1$ and $x_2$ in $\R_+$, $f(x_1,
g(x_1)) = 0$ and $f(k(x_2), x_2) = 0$, it is clear that $W
\subset Z$. Let $q := (q_1, q_2) \in Z$.

Suppose that $q_1 < p_1$ and $q_2 < p_2.$ Let
$$z_1 := \frac{\log(q_2/p_2)}{\log(p_1/q_1)}$$
So, $p_1^{z_1}p_2 = q_1^{z_1}q_2$. As $p_1/q_1 > 1$ and $q_2 /
p_2 < 1$, then $z_1 < 0$. Let $z \in \R^2, z := (z_1, 1)$. We
know that $p_1$ is a zero of multiplicity at least $2$ of
$f_{(p,z)}$. On the other hand,
$$f_{(p,z)}(q_1) = f(q_1, p^zq_1^{-z_1}) = f(q_1, q^zq_1^{-z_1}) = f(q_1, q_2) =  0,$$
because $q \in Z$. Then, $f_{(p,z)}$ has at least three zeros
(counted with multiplicity) and, by \DR, at least three sign
changes. As $c, d > 1$ and $z_1 < 0$, then $0 <  -z_1 < c - dz_1
 \hbox{ and } 0 < 1 <  c - dz_1$, and
we have that
$$f_{(p,z)}(x_1)  =  1 - x_1 - p^zx_1^{-z_1} + A(p^z)^{d}x_1^{c - dz_1}$$
has just two sign changes. Then, it cannot happen that $q_1 <
p_1$ and $q_2 < p_2$ at the same time.

Suppose now that $q_1 \ge p_1$. Consider the following cases:
\begin{itemize}
\item $q_1 \ge 1$: as we have shown at the beginning of this
Lemma, the line $\{x_1 = q_1\}$ intersects $Z$ in a single point,
which is $(q_1, g(q_1))$. Then, it must be $q_2 = g(q_1)$ and then
$q \in W_1$. \item $p_1 <q_1 < 1$: we know that the line $\{x_1 =
q_1\}$ intersects $Z$ in two points: $(q_1, g(q_1))$ and $(q_1,
t(q_1))$, with $g(q_1) \in ((p^v)^{1/(d-1)}x_1^{-c/(d-1)}, p_2)$
and $t(q_1) \in (0, (p^v)^{1/(d-1)}q_1^{-c/(d-1)})$. If $q_2 =
g(q_1)$ then $q \in W_1$. If $q_2 = t(q_1)$, then
$$q_2 < (p^v)^{1/(d-1)}q_1^{-c/(d-1)} =
(p_1/q_1)^{c/d-1}p_2 < p_2,$$ and therefore $x_1 = k(q_2)$ is the
unique value of $x_1$ in the interval $(p_1,$
$(p^w)^{1/(c-1)}q_2^{-d/(c-1)})$ such that $f(x_1, q_2) = 0$.
Since the previous inequalities imply
$$q_1 < p_1p_2^{(d-1)/c}q_2^{-(d-1)/c}
< p_1p_2^{d/(c-1)}q_2^{-d/(c-1)} = (p^w)^{1/(c-1)}q_2^{-d/c-1},$$
we conclude that $q_1 = k(q_2)$ and so $q \in W_2$. \item  If
$q_1 = p_1$, let us see that $\{x_1 = q_1\}$ intersects $Z$ only
in $p$. Let us consider $e_1 = (1,0)$. We know that  $p_2$ is a
zero of multiplicity at least $2$ of  $f_{(p, e_1)}$, but
$$f_{(p, e_1)}(x_2) =  (1-p_1) - x_2 + Ap_1^cx_2^d$$
is a $3$-nomial with two sign changes. By \DR , $p_2$ has
multiplicity equal to $2$ and $f_{(p, e_1)}$ has no other zeros.
Then $\{x_1 = q_1\} \cap
Z = \{p\}$, and then $q = p \in Z$.
\end{itemize}

If $q_2 \ge p_2$, we proceed in an analogous way. Thus, we
conclude that $Z = W$, and that $Z$ has a unique connected
component, which is unbounded.
\end{proof}

We can now give a proof of the following theorem.

\begin{thm}
\label{final} Let $f$ be a $4$-nomial in two variables and let $Z
:= f^{-1}(0) \subset \R^2_+$. If $Z$ has a compact connected
component $\Gamma$, then $Z = \Gamma$.
\end{thm}

\begin{proof}{Proof:}

Suppose $Z \setminus \Gamma \ne \emptyset$. By Lemma
\ref{SobranCosas}, we know that among the coefficients of $f$
there are two positive and two negative. On the other hand, we
know that $\dim \mathrm{Newt}(f) = 2$, otherwise $Z$ could not
have compact connected components. Then, by Lemma \ref{CuadriOpu},
$\mathrm{Newt}(f)$ is a quadrilateral without parallel sides and
coefficients of the same sign correspond to opposite vertices. By
Lemma \ref{PtoCrit}, $Z$ does not have critical points; otherwise
it would have only a unique non-compact connected component.
Finally, because of Lemma \ref{Equiv}, we can suppose $f$ is of
the following type:
$$f(x_1, x_2) = 1 - x_1 - x_2 + Ax_1^cx_2^d,$$
with $A > 0; c, d >1.$

Again, in order to study how many times the line $\{x_1 =
\alpha\}$ intersects $Z$ for a fixed $\alpha$ in $\R_+$, let us
define a function $g_{\alpha}$ in the variable $x_2$ as the
restriction of $f$ to that line, i.e.:
$$g_{\alpha}(x_2) = f(\alpha, x_2).$$
Then
$$g'_{\alpha}(x_2) = - 1 + Ad\alpha^c x_2^{d-1} < 0
\iff x_2 < \Big(\frac{1}{Ad}\Big)^{1/(d-1)}\alpha ^{-c/(d-1)}.$$
Let $J := (1/Ad)^{1/(d-1)}$. Then $J > 0$ and the function
$g_{\alpha}$ has a minimum in $x_2 = J\alpha^{-c/(d-1)}$. For $x_1
\in \R_+$, let $\ell_1 (x_1)$ be the minimum of the function
$g_{x_1}$ and let $\ell (x_1):= (x_1, \ell_1 (x_1))$. Then,
$$f \circ \ell (x_1) = (- J + AJ^d)x_1^{-c/(d-1)} + 1 - x_1,$$
so $f \circ \ell $ turns out to be a $3$-nomial.

As $\Gamma$ is a compact set, the function $x_1$ reaches its
minimum (let us call it $m$) and its maximum (let us call it $M$)
on $\Gamma$. Let us prove that $m\ne M$: as $\Gamma$ is a
differentiable manifold of dimension $1$, then $\Gamma$ has an
infinite number of points. If $m = M$, then $\Gamma \subset \{x_1
= m \}$, and the $3$-nomial $g_m$ has infinitely many zeros,
which is impossible.

Let $p := (p_1, p_2)$ and $q:= (q_1, q_2)$ in $\Gamma$ such that
$m = p_1$ and $M = q_1$. Let us see that $p_2$ is a zero of
multiplicity $2$ of $g_{p_1}$. As $p \in \Gamma$, $g_{p_1}(p_2) =
f(p_1, p_2) = 0$. On the other hand, as $p$ is the minimum of
$x_1$ in $\Gamma$, using Lagrange multipliers,
$$g'_{p_1}(p_2) = \frac{\partial f}{\partial x_2}(p) = 0.$$
By \DR, we know that $g_{p_1}$ must have at least two sign
changes. As we know that
$$g_{p_1}(x_2) = (1 - p_1) - x_2 + Ap_1^cx_2^d,$$
then it must be $p_1 < 1$, and so $g_{p_1}$ has no zeros other
than $p_2$. As the unique zero of $g_{p_1}$ has an even
multiplicity and its leading coefficient is positive, for all
$x_2 \ne p_2$, $g_{p_1}(x_2)
> 0$. Then, $f \circ \ell(p_1) = 0$. In an analogous way, we can
prove that $q_2$ is the unique zero of the function $g_{q_1}$ and
$f \circ \ell (q_1) = 0$. We conclude that $p_1 = m$ and $q_1 = M$
are two different zeros of $f \circ \ell$.

As $f\circ \ell (x_1) = (- J + AJ^d)x_1^{-c/(d-1)} + 1 - x_1$ is
a $3$-nomial, 
it has no zeros other than $m$ and $M$ which have multiplicity
$1$. As its leading coefficient is negative, we know that $f \circ
\ell(x_1) < 0$ for all $x_1 \in (0,m) \cup (M, +\infty)$ and $f
\circ \ell (x_1) > 0$ for all $ x_1 \in (m, M)$. Let $s \in (m,
M)$. Then, for all $x_2 \in \R_+$, $f(s, x_2) \ge f \circ \ell(s)
> 0$,  . Then, $\Gamma \cap \{x_1 = s\} = \emptyset$ and the open
sets $\{x_1 < s\}$ and $\{x_1 > s\}$ disconnect $\Gamma$, which
is a contradiction.
\end{proof}

Now we can give a proof of Theorem \ref{Resfinal}, which is the
main goal of this section.

\begin{proof}{Proof:}
$\hbox{

}$
\begin{enumerate}
\item The inequality $P_{comp}(2,4) \le 1$ is a consequence of
Theorem \ref{final}. In the following example the equality holds:
$$f_1(x_1, x_2) = x_2^2 - 4x_1^3x_2 + x_1^8 + 3x_1^4.$$
In fact, $f_1(x_1, x_2) = 0$ if and only if $x_2 = 2x_1^3 \pm
x_1^2 \sqrt{1 - (x_1^2 - 2)^2}$, and the set of positive values
of $x_1$ where the polynomial under the square root symbol is
non-negative is the interval $[1 , \sqrt3]$.

\item Let $f$ be a $4$-nomial in two variables and let $Z :=
f^{-1}(0) \subset \R^2_+.$

If $\dim \mathrm{Newt}(f)= 1$, then by Proposition \ref{dimcorr}
and \DR, we know that $\mathrm{Non}(Z) \le P(1,4) \le 3$.

If $\dim \mathrm{Newt}(f)= 2$ and $0$ is a regular value of $f$,
then by \cite[Theorem 3]{LiRojWang}, $\mathrm{Non}(Z) \le 2$.

If $\dim \mathrm{Newt}(f)= 2$ and $0$ is not a regular value of
$f$, there is
 a critical point $p$ in $Z$. If $Z = \{p\}$, then $\mathrm{Non}(Z) = 0$.
If $Z \ne \{p\}$, by Lemma \ref{PtoCrit}, $\mathrm{Non}(Z) = 1$.

The equality holds in the following example:
$$f_2(x_1, x_2) = (x_1 - 1)(x_1 - 2)(x_1 - 3) = x_1^3 - 6x_1^2 + 11x_1 - 6.$$

\item Let $f$ be a $4$-nomial in two variables, and let $Z :=
f^{-1}(0) \subset \R^2_+.$

If $Z$ has any compact connected component $\Gamma$, by Theorem
\ref{final}, $Z = \Gamma$ and then $\mathrm{Tot}(Z) = 1$. If it
does not, because of the previous item we have that
$\mathrm{Tot}(Z) \le 3$ and the same example shows that the
equality holds.

\item Let $f$ be a $4$-nomial in two variables such that $\dim
\mathrm{Newt}(f) = 2$, and let $Z := f^{-1}(0) \subset \R^2_+.$

If $Z$ has any compact connected component $\Gamma$, again by
Theorem \ref{final}, $Z = \Gamma$ and then $\mathrm{Tot}(Z) = 1$.
If it does not, as it was shown in the second item of this
theorem, $\mathrm{Non}(Z) \le 2$, and $\mathrm{Tot}(Z) \le 2$.

The equality holds in the following example:
$$f_3(x_1, x_2) =
x_1x_2 - 2x_1 - x_2 + 1.$$ In fact, $f(x_1, x_2) = 0$ is an
implicit equation for the hyperbola $x_2 = \frac1{x_1 - 1} + 2$.
\end{enumerate}
\end{proof}

\section{On $m$-nomials in $n$ variables}

In this section we will prove Theorems \ref{CuentaFinal} and
\ref{Anexo}. Theorem \ref{CuentaFinal} gives us an explicit upper
bound for the number of connected components of the zero set of
an $m$-nomial in $n$ variables in the positive orthant and Theorem
\ref{Anexo} is an auxiliary theorem for Theorem
\ref{CuentaFinal}, but it also will be used in the next section.
Let us give now a proof of Theorem \ref{Anexo}.

\begin{proof}{Proof:}

As observed earlier in Remark \ref{sepuededividir}, we can assume
$f$ is of the following form:
$$f(x) = c_m + \sum_{i=1}^{m-1}c_ix^{a_i}.$$
As $\dim\langle a_1$, $\dots$, $a_{m-1}\rangle = n$, without loss
of generality, we can suppose that $B : =  \{a_1, \dots, a_n\}$
is a basis of $\R^n$. Let $A$ be the matrix having the elements
of $B$ as columns, let $g$ be the $m$-nomial $f \circ h_{A^{-1}}$
and let $W := g^{-1}(0) \subset \R^n_+$. As $h_{A^{-1}}$ is a
diffeomorphism,  $\mathrm{Non}(W) = \mathrm{Non}(Z)$ and $\dim
\mathrm{Newt}(g) = \dim \mathrm{Newt}(f)$. Moreover, we have that
$$g(x) = c_m + \sum_{i=1}^{m-1} c_ix^{A^{-1}a_i}
 = c_m + \sum_{i=1}^{n} c_ix_i + \sum_{i=n+1}^{m-1} c_ix^{A^{-1}
a_i}.$$

Suppose $W$ has $t$ non-compact connected components and let $\{
p_1, \dots, p_t\}$ be a set of points intersecting each and every
non-compact connected component of $W$. Suppose, for $1 \le i \le
t$, $p_i = (p_{i1}, \dots, p_{in})$. For $1\le j \le n$, let us
consider $M_j, m_j \in \R_+$ such that $M_j > \max\{ p_{ij}, \,1
\le i \le t\}$ and $m_j < \min\{ p_{ij}, \,1 \le i \le t\}$, and
define $S_j = \{x \in \R_+^n \; | \; x_j = M_j \}$ and $T_j = \{x
\in \R_+^n \; | \; x_j = m_j \}$. Let us prove that each
non-compact connected component of $W$ intersects at least one of
the sets $S_1, \dots, S_n, T_1, \dots, T_n$.

Let $X$ be a non-compact connected component of $W$. If $X$ is
not bounded, then there exists $j_0, 1 \le j_0 \le n,$ such that
$X \cap S_{j_0}$ is not empty.

If $X$ is bounded, then it is not closed. Let $T := \cap_{j = 1}^n
\{x \in \R_+^n \; | \; x_j \ge m_j \}.$ If $X \subseteq T$, then
it is a connected component of $W \cap T$. As $W = g^{-1}(0)
\subset \R^n_+$ and $g$ is a continuous function, there exists a
closed set $F \subset \R^n $ such that $W = F \cap \R^n_+$. Then
$$W \cap T =F \cap \R^n_+ \cap T = F \cap T,$$
and $W \cap T \subset \R^n_+ $ is closed because it is an
intersection of closed sets. It follows that $X$ is closed because
it is a connected component of a closed set. This is a
contradiction, and then $X \nsubseteq T$, and this implies that
there exists $j_1$, $1 \le j_1 \le n,$ such that $X \cap T_{j_1}
\ne \emptyset$.

In this way, we have found $2n$ sets ($S_1, \dots, S_n, T_1,
\dots, T_n)$ such that each non-compact connected component of
$W$ has a non-empty intersection with one of them. Thus,
$$\mathrm{Non}(W)  \le \sum_{j = 1}^n
 \mathrm{Tot}(W \cap S_j)   +  \sum_{j = 1}^n  \mathrm{Tot}(W \cap T_j).$$

Each of these $2n$ intersections has at most $P(n-1, m-1)$
connected components, because they can be regarded as zero sets of
$m'$-nomials in $n-1$ variables, with $1 \le m' \le m-1$. For
example, the set $W \cap S_n$ can be described as the zero set of
the following function:
$$\hat g: \R_+^{n-1} \to \R,$$
$$\hat g(x_1, \dots, x_{n-1})
= (c_m + c_nM_n)+ \sum_{i=1}^{n-1} c_ix_i + \sum_{i=n+1}^{m-1}
c_i(x_1, \dots, x_{n-1},M_n)^{A^{-1}a_i}.$$

We have thus proved that
$$\mathrm{Non}(Z) = \mathrm{Non}(W)
\le 2n\,P(n-1, m-1),$$ which is our first assertion.

Finally, note that for the function $\hat g$ defined above, $\dim
\mathrm{Newt}(\hat g) = n-1$. Proceeding inductively, we get the
second inequality.
\end{proof}

Let us prove now theorem \ref{CuentaFinal}.

\begin{proof}{Proof:}

Let us proceed by induction on $n$.

If $n = 1$ by \DR, we know that
$$P(1,m) \le m-1 < 2^{m-1}2^{1 + (m-1)(m-2)/2}.$$

Suppose now that $n > 1$. Given an $m$-nomial $f$ in $n$
variables, let $d := \dim \mathrm{Newt}(f)$ and $Z := f^{-1}(0)
\subset \R^n_+$. 

If $d < n$, by the first item of Proposition \ref{dimcorr} and the
induction hypothesis,
$$\mathrm{Tot}(Z) \le P(d,m)
\le (d+1)^{m-1}2^{1 + (m-1)(m-2)/2} \le (n+1)^{m-1}2^{1 +
(m-1)(m-2)/2}.$$

If $d = n$, as $m-1 \ge d$, we have that $m \ge n + 1$. If $m = n
+ 1$, by the second item of Proposition \ref{dimcorr},
$\mathrm{Tot}(Z) \le 1$. If $m \ge n + 2$, by the second item of
Theorem \ref{Anexo}, 
the first item of Theorem \ref{TeoLiRojWang} and Theorem
\ref{Khovanski}, we have
$$\mathrm{Tot}(Z) \le \sum_{i = 0}^{n-1}
2^{i} \frac{n!}{(n - i)!}(n-i+1)^{m-i-1}2^{(m-i-1)(m-i-2)/2}.$$

Now, we use the following inequality, valid for all $i, n, m \in
\N$ such that $m \ge n+2$ and $0 \le i \le n-1$, that can be
easily proved by induction on $i$:
$$2^i\frac{n!}{(n-i)!}(n-i+1)^{m-i-1}2^{(m-i-1)(m-i-2)/2}
\le \frac{1}{2^i}(n+1)^{m-1}2^{(m-1)(m-2)/2}.$$ Then, we conclude
that
$$\mathrm{Tot}(Z) \le \sum_{i = 0}^{n-1}
\frac{1}{2^i}(n+1)^{m-1}2^{(m-1)(m-2)/2} <  (n+1)^{m-1}2^{1 +
(m-1)(m-2)/2},$$ which completes the proof.
\end{proof}

\section{On $5$-nomials in three variables}

As a consequence of what has been proved in the previous sections,
we get Theorem \ref{five}:

\begin{proof}{Proof:}
By Theorem \ref{Anexo}, $\mathrm{Non}(Z) \le 6P(2,4) = 18.$
Nevertheless, in the proof of that theorem, we have shown the
existence of six $4$-nomials in two variables, let us call them
$g_1, \dots, g_6$, such that for $i = 1, \dots, 6$, $\dim
\mathrm{Newt}(g_i) = 2$ and
$$\mathrm{Non}(Z) \le \sum_{i = 1}^6 \mathrm{Tot}(g_i^{-1}(0)).$$
By the fourth item of Theorem \ref{Resfinal}, we know that
$\mathrm{Tot}(g_i^{-1}(0)) \le 2$, and then we conclude that
$\mathrm{Non}(Z) \le 12.$
\end{proof}

This bound is significantly sharper than the best previously known
one, which was $10384$.

\bigskip

\noindent {\bf Acknowledgments:} The author thanks Teresa Krick
and Juan Sabia for their guidance in this work, Fernando L\'opez
Garc\'\i a for his valuable cooperation, Maurice Rojas for some
nice talks and Gabriela Jeronimo and the referees for their very
helpful corrections and suggestions.

\end{document}